\documentclass[a4paper,12pt]{article}
\usepackage[top=2.5cm,bottom=2.5cm,left=2.5cm,right=2.5cm]{geometry}
\usepackage{cite, amsmath, amssymb}
\usepackage[margin=1cm,%
            font=small,%
            format=hang,%
            labelsep=period,%
            labelfont=bf]{caption}
\numberwithin{equation}{section}
\usepackage{graphicx}
\usepackage{subfigure}
\usepackage{cite}
\allowdisplaybreaks[1]

\begin{document}
\title{\vspace{3.5 cm} The (multiplicative degree-)Kirchhoff index of graphs derived from the Catersian product of $S_n$ and $K_2$
\footnote{Supported by the NSFC (11771362) and the Natural Science Foundation of Anhui Province (2008085J01)}
}

\author{Jia-Bao Liu $^1$, Xin-Bei Peng $^1$, Jiao-Jiao Gu $^1$,
Wenshui Lin $^{2,3,}$\footnote{Corresponding author. Email address: wslin@xmu.edu.cn (W. Lin).}\vspace*{0.8cm}\\
    {\small $^1$ School of Mathematics and Physics, Anhui Jianzhu University, Hefei 230601, China}\\
    {\small $^2$ Fujian Key Laboratory of Sensing and Computing for Smart City, Xiamen 361005, China}\\
    {\small $^3$ School of Informatics, Xiamen University, Xiamen 361005, China}
}

\date{\small (Received July 21, 2020)}
\maketitle \thispagestyle{empty}

\noindent \textbf{Abstract}\\
Recently, Li et al. [Appl. Math. Comput. 382 (2020) 125335] proposed the problem
of determining the Kirchhoff index and multiplicative degree-Kirchhoff index of graphs derived from $S_n \times K_2$,
the Catersian product of the star $S_n$ and the complete graph $K_2$.
In the present paper, we completely solve this problem. That is,
the explicit closed-form formulae of Kirchhoff index, multiplicative degree-Kirchhoff index,
and number of spanning trees are obtained for some graphs derived from $S_n \times K_2$.\\


\baselineskip=0.30in

\section{Introduction}
Let $G=(V,E)$ be a nontrivial simple connected graph,
where $V = \{v_1, v_2, \cdots, v_n\}$ and $E$ are the vertex set and edge set of $G$, respectively.
$A(G)=(a_{ij})_{n \times n}$ is the adjacency matrix of $G$,
where $a_{ij} =1$ if $v_i v_j \in E$, and $a_{ij} =0$ otherwise.
Let $d_i$ be the degree of vertex $v_i$ in $G$, and $D(G) = diag(d_1, d_2, \cdots, d_n)$.
Then $L(G) =D(G) -A(G)$ is called the Laplacian matrix of $G$,
and $\mathcal{L}(G) =D(G)^{-\frac{1}{2}} L(G) D(G)^{-\frac{1}{2}}$ the normalized Laplacian matrix of $G$.
It is easily seen that,
\begin{equation*}
  \left(\mathcal{L}(G)\right)_{ij} = \begin{cases}
    1,&\text{if}~i =j;\\
    -\frac{1}{\sqrt {d_i d_j}}, &\text{if}~ i \ne j ~\text{and}~ v_i v_j \in E;\\
    0, &\text{otherwise}.
   \end{cases}
\end{equation*}

Let $0= \mu_1 <\mu_2 \leq \cdots \leq \mu_n$ be the eigenvalues of $L(G)$,
and $0= \nu_1 <\nu_2 \leq \cdots \leq \nu_n$ the eigenvalues of $\mathcal{L}(G)$.
The sets $Sp(L(G)) = \{\mu_1, \mu_2, \cdots, \mu_n\}$
and $Sp(\mathcal{L}(G)) = \{\nu_1, \nu_2, \cdots, \nu_n\}$ are called
the Laplacian spectrum and normalized Laplacian spectrum of $G$, respectively.

Let $d_{ij}$ denote the distance between vertices $v_i$ and $v_j$ in $G$
(namely, the length of a shortest path connecting them).
The Wiener index \cite{b1} and Gutman index \cite{b2} of $G$ are defined as
$W(G) =\sum_{i<j} d_{ij}$ and $Gut(G) =\sum_{i<j} d_i d_j d_{ij}$.
For these two famous topological indices, one can refer to \cite{b3,b4,b5,b6,b7,b8,b9,b10}
and the references therein.

If regard each edge in $E(G)$ as an unit resistor,
then the resistance distance between two vertices $v_i$ and $v_j$,
denoted by $r_{ij}$, is defined \cite{b11} to be the effective resistance between them.
Similar to the Wiener index, the Kirchhoff index of $G$ is defined as $Kf(G)=\sum_{i<j} r_{ij}$.
Later, the following relation between $Kf(G)$ and $Sp(L(G))$ was established
by Zhu et al. \cite{b12} and Gutman and Mohar \cite{b13} independently.

\noindent \textbf{Lemma 1.1} \cite{b12,b13}\textbf{.} Let $G$ be a simple graph of order $n \geq 2$.
Then
$$Kf(G) = \sum_{i=2}^n \frac{1}{\mu_i}.$$

Similar to the Gutman index, Chen and Zhang \cite{b14} defined
the multiplicative degree-Kirchhoff index of $G$ as
$Kf^*(G)= \sum_{i <j} d_i d_j r_{ij}$.
Moreover, the following relation between $Kf^*(G)$ and $Sp(\mathcal{L}(G))$ were confirmed.

\noindent \textbf{Lemma 1.2} \cite{b14}\textbf{.}
Let $G$ be a simple connected graph of order $n \geq 2$ and size $m$. Then
$$Kf^*(G) = 2m \sum_{i=2}^n \frac{1}{\nu_i}.$$

In recent years, more and more attentions were paid to the Kirchhoff index and
multiplicative degree-Kirchhoff index.
Closed-form formulae of (multiplicative degree-)Kirchhoff index
have been established for some classes of graphs.
For examples, the formulae of Kirchhoff index for cycles, circulant graphs,
and composite graphs were obtained in \cite{b15}, \cite{b16}, and \cite{b17}, respectively,
and those of both indices for complete multipartite graphs were obtained in \cite{b18}.
Besides, quite a few literatures concerned the (multiplicative degree-)Kirchhoff index
of polygon chains and their variants.
Explicit expressions of (multiplicative degree-)Kirchhoff index
have been obtained for linear polyomino chain \cite{b19},
linear crossed polyomino chain \cite{b20}, linear pentagonal chain \cite{b21},
linear phenylenes \cite{b22,b23}, cyclic phenylenes \cite{b24},
M\"{o}bius phenylenes chain and cylinder phenylenes chain \cite{b25,b26},
linear [n] phenylenes \cite{b27}, generalized phenylenes \cite{b28,b29},
linear hexagonal chain \cite{b30,b31}, linear crossed hexagonal chain \cite{b32},
M\"{o}bius hexagonal chain \cite{b33}, and periodic linear chains \cite{b34},
linear octagonal chain \cite{b35}, linear octagonal-quadrilateral chain \cite{b36},
and linear crossed octagonal chain \cite{b37}.

For two disjoint graphs $G$ and $H$, $G \otimes H$ will denote the strong product of $G$ and $H$.
That is, $V(G \otimes H) = V(G) \times V(H)$,
and two distinct vertices $(u_1, v_1)$ and $(u_2, v_2)$ are adjacent
whenever $u_1 =u_2$ or $u_1 u_2 \in E(G)$, or, $v_1 =v_2$ or $v_1 v_2 \in E(H)$.
The Catersian product of $G$ and $H$, denoted by $G \times H$,
is the graph with vertex set $V(G) \times V(H)$,
and two vertices $(u_1, v_1)$ and $(u_2, v_2)$ are adjacent
whenever $u_1= u_2$ and $v_1 v_2 \in E(H)$, or $v_1= v_2$ and $u_1 u_2 \in E(G)$.
Figure 1 depicts the graphs $S_n \otimes K_2$ and $S_n \times K_2$,
where $S_n$ and $K_n$ denote the star and complete graph of order $n$, respectively.
Recently, Li et al. \cite{b38} determined the expressions of $Kf(S_r)$,
$Kf^*(S_r)$, and $\tau(S_r)$, where $S_r$ is a graph derived from $S_n \otimes K_2$
by randomly removing $r$ vertical edges,
and $\tau(G)$ denotes the number of spanning trees of a connected graph $G$.
Finally, they proposed the problem of determining these three invariants for
graphs derived from $S_n \times K_2$. In the present paper, we completely solve this problem.

\begin{figure}[ht]
 \centering
 \includegraphics[width=4 in]{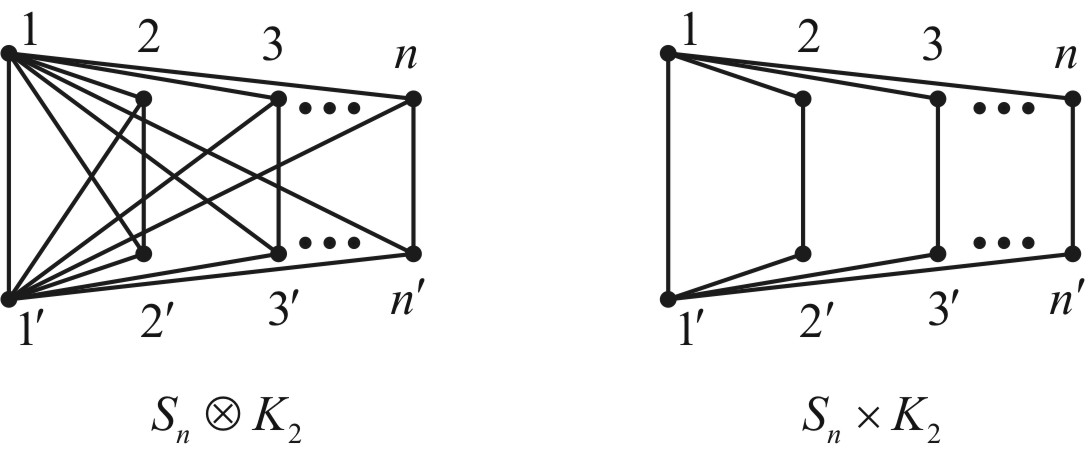}
 \caption{ The graphs $S_{n} \otimes K_{2}$ and $S_{n} \times K_{2}$.}
\end{figure}

For convenience, we denote $S_n^2 = S_n \times K_2$.
Then $|V(S_n^2)| = 2n$ and $|E(S_n^2)| = 3n-2$.
Let $E'= \{ii'|i=1,2,\cdots,n\}$. 
$\mathcal{S}_{n,r}^2$ will denote the set of graphs derived from $S_n^2$ by
randomly deleting $r$ edges in $E'$.
Obviously, the unique graph in $\mathcal{S}_{n,n}^2$ is disconnected,
hence we consider $\mathcal{S}_{n,r}^2$ for $0 \leq r \leq n-1$ only.
Note also, $\mathcal{S}_{n,0}^2 = \{ S_n^2 \}$.
In Section 2, some notations and known results are introduced,
which will be applied to get our main results.
In Section 3, explicit expressions of $Kf(S_n^2)$, $Kf^*(S_n^2)$, and $\tau(S_n^2)$ are obtained.
Finally, $Kf(S_{n,r}^2)$ and $\tau(S_{n,r}^2)$ are determined in Section 4,
where $S_{n,r}^2$ is an arbitrary graph in $\mathcal{S}_{n,r}^2$.
Moreover, it is shown that, $\lim_{n \rightarrow +\infty} Kf(S_n^2)/W(S_n^2) = \lim_{n \rightarrow +\infty} Kf(S_{n,r}^2)/W(S_{n,r}^2) = 8/15$
and $\lim_{n \rightarrow +\infty} Kf^*(S_n^2)/Gut(S_n^2) = 16/33$.

\section{Preliminaries}
Label the vertices of $S_n^2$ as in Figure 1, and set $V_1= \{1,2,\cdots,n\}$,
$V_2= \{1',2',\cdots,n'\}$. Then we have
$$L(S_n^2)= \begin{pmatrix} L_{11}(S_n^2) & L_{12}(S_n^2)\\
                            L_{21}(S_n^2) & L_{22}(S_n^2)
            \end{pmatrix},~
  \mathcal{L}(S_n^2)= \begin{pmatrix} \mathcal{L}_{11}(S_n^2) & \mathcal{L}_{12}(S_n^2)\\
                                      \mathcal{L}_{21}(S_n^2) & \mathcal{L}_{22}(S_n^2)
            \end{pmatrix},$$
\noindent where $L_{ij}(S_n^2)$ ($\mathcal{L}_{ij}(S_n^2)$)
is the submatrix of $L(S_n^2)$ (resp. $\mathcal{L}(S_n^2)$)
whose rows (columns) corresponding to the vertices in $V_i$ (resp. $V_j$).
It is easily seen that, $L_{11}(S_n^2) =L_{22}(S_n^2)$,
$L_{12}(S_n^2)= L_{21}(S_n^2)$, $\mathcal{L}_{11}(S_n^2) = \mathcal{L}_{22}(S_n^2)$,
and $\mathcal{L}_{12}(S_n^2) = \mathcal{L}_{21}(S_n^2)$.

Let
$$ T=\begin{pmatrix}
     \frac{1}{\sqrt{2}}I_n &\frac{1}{\sqrt{2}}I_n \\
     \frac{1}{\sqrt{2}}I_n &-\frac{1}{\sqrt{2}}I_n
     \end{pmatrix},$$

\noindent then we have
$$T L(S_n^2) T = \begin{pmatrix}
                    L_A(S_n^2) &0\\
                    0 &L_S(S_n^2)
                 \end{pmatrix},~
 T \mathcal{L}(S_n^2) T = \begin{pmatrix}
                    \mathcal{L}_A(S_n^2) &0\\
                    0 &\mathcal{L}_S(S_n^2)
                 \end{pmatrix},$$

\noindent where
$L_A(S_n^2)= L_{11}(S_n^2) + L_{12}(S_n^2)$, $L_S(S_n^2)= L_{11}(S_n^2) - L_{12}(S_n^2)$,
$\mathcal{L}_A(S_n^2) = \mathcal{L}_{11}(S_n^2)+ \mathcal{L}_{12}(S_n^2)$,
and $\mathcal{L}_S(S_n^2) = \mathcal{L}_{11}(S_n^2) - \mathcal{L}_{12}(S_n^2)$.

Based on the above arguments, by applying the technique used in \cite{b32,b39},
we immediately have the following decomposition theorem,
where $\Phi(B, \lambda) = |\lambda I-B|$ stands for
the characteristic polynomial of $B$.

\noindent \textbf{Lemma 2.1.} Let $L_A(S_n^2)$, $L_S(S_n^2)$, $\mathcal{L}_A(S_n^2)$,
and $\mathcal{L}_S(S_n^2)$ be defined as above.
Then
$$\Phi(L(S_n^2), \lambda) =  \Phi(L_A(S_n^2), \lambda) \Phi(L_S(S_n^2), \lambda),$$
and
$$\Phi(\mathcal{L}(S_n^2), \lambda)
= \Phi(\mathcal{L}_A(S_n^2), \lambda) \Phi(\mathcal{L}_S(S_n^2), \lambda).$$

\noindent \textbf{Lemma 2.2} \cite{b40}\textbf{.} If $G$ is a connected graph with $n \geq 2$ vertices,
then $$\tau(G) = \frac{1}{n}\prod_{i=2}^{n} \mu_i.$$

\section{Results for $S_n^2$}
We will give explicit expressions of $Kf(S_n^2)$, $Kf^*(S_n^2)$, and $\tau(S_n^2)$ in this section.

\subsection{On $Kf(S_n^2)$ and $\tau(S_n^2)$}
Obviously,
$$L_{11}(S_n^2) = \begin{pmatrix}
                   n &-1 &-1 &\cdots &-1 \\
                   -1 &2 &0 &\cdots &0 \\
                   -1 &0 &2 &\cdots &0\\
                   \cdots &\cdots &\cdots &\cdots &\cdots\\
                   -1 &0 &0 &\cdots &2
                  \end{pmatrix}_{n \times n},~
  L_{12}(S_n^2) = \begin{pmatrix}
                   -1 &0 &0 &\cdots &0 \\
                   0 &-1 &0 &\cdots &0 \\
                   0 &0 &-1 &\cdots &0\\
                   \cdots &\cdots &\cdots &\cdots &\cdots\\
                   0 &0 &0 &\cdots &-1
                  \end{pmatrix}_{n \times n}.$$

Hence
$$L_A(S_n^2) = L_{11}(S_n^2)+ L_{12}(S_n^2) = \begin{pmatrix}
                   n-1 &-1 &-1 &\cdots &-1 \\
                   -1  &1  &0  &\cdots &0 \\
                   -1  &0  &1  &\cdots &0\\
                   \cdots &\cdots &\cdots &\cdots &\cdots\\
                   -1  &0  &0  &\cdots &1
                  \end{pmatrix}_{n \times n},$$
\noindent and we easily have $Sp(L_A(S_n^2)) = \{0, 1^{n-2}, n \}$,
where $a^k$ denotes $k$ successive $a$'s.

Similarly, we have
$$L_S(S_n^2) = L_{11}(S_n^2) - L_{12}(S_n^2) = \begin{pmatrix}
                   n+1 &-1 &-1 &\cdots &-1 \\
                   -1  &3  &0  &\cdots &0 \\
                   -1  &0  &3  &\cdots &0\\
                   \cdots &\cdots &\cdots &\cdots &\cdots\\
                   -1  &0  &0  &\cdots &3
                  \end{pmatrix}_{n \times n},$$
\noindent and get $Sp(L_S(S_n^2)) = \{2, 3^{n-2}, n+2 \}$.

Hence $Sp(L(S_n^2)) = \{0, 1^{n-2}, 2, 3^{n-2}, n, n+2 \}$ from Lemma 2.1,
and we get the following result.

\noindent \textbf{Theorem 3.1.} Let $S_n^2 = S_n \times K_2$. Then

(1) $Kf(S_n^2) = \frac{8n^3 +3n^2 -14n +12}{3n +6}$;

(2) $\tau(S_n^2) = (n+2)\cdot 3^{n-2}$;

(3) $\lim\limits_{n \rightarrow +\infty} \frac{Kf(S_n^2)}{W(S_n^2)} = \frac{8}{15}$.

\noindent \textbf{Proof.} From Lemma 1.1 we have
$$Kf(S_n^2) = 2n \left[(n-2) + \frac{1}{2} + \frac{n-2}{3} + \frac{1}{n} + \frac{1}{n+2} \right]
= \frac{8n^3 +3n^2 -14n +12}{3(n+2)}.$$

From Lemma 2.2 we immediately have
$$\tau(S_n^2) = \frac{1}{2n}\cdot 2 \cdot 3^{n-2} \cdot n \cdot (n+2) = (n+2)\cdot 3^{n-2}.$$

Finally, we end the proof by confirm that $W(S_n^2) =5n^2 -8n +4$.
Let $w_i = \sum_{j \in V(S_n^2)} d_{ij}$.
Obviously, $w_i = 1 \cdot n + 2(n-1) = 3n-2 $  if  $i = 1, 1'$,
and $w_i = 1 +1 +2(n-1) +3(n-2) = 5n -6$ otherwise.
Hence $$ W(S_n^2) = \frac{1}{2} \sum\limits_{i \in V(S_n^2)} w_i
                  = \frac{1}{2} \left[ 2(3n-2) +(2n-2)(5n-6) \right] = 5n^2 -8n +4.~\blacksquare$$

\subsection{On $Kf^*(S_n^2)$}
Consequently, we will determine $Kf^*(S_n^2)$. Obviously,
$$\mathcal{L}_{11}(S_n^2) =\begin{pmatrix}
        1 &-\frac{1}{\sqrt{2n}} &-\frac{1}{\sqrt{2n}} &\cdots &-\frac{1}{\sqrt{2n}}\\
        -\frac{1}{\sqrt{2n}} &1 &0 &\cdots &0\\
        -\frac{1}{\sqrt{2n}} &0 &1 &\cdots &0\\
        \cdots &\cdots &\cdots &\cdots &\cdots\\
        -\frac{1}{\sqrt{2n}} &0 &0 &\cdots &1
   \end{pmatrix}_{n \times n},
$$
\noindent and
$$\mathcal{L}_{12}(S_n^2) =\begin{pmatrix}
        -\frac{1}{n} &0 &0 &\cdots &0 \\
        0 &-\frac{1}{2} &0 &\cdots &0\\
        0 &0 &-\frac{1}{2} &\cdots &0\\
        \cdots &\cdots &\cdots &\cdots &\cdots\\
        0 &0 &0 &\cdots &-\frac{1}{2}
   \end{pmatrix}_{n \times n}.
$$
\noindent Hence
$$\mathcal{L}_A(S_n^2) = \mathcal{L}_{11}(S_n^2) + \mathcal{L}_{12}(S_n^2) =
   \begin{pmatrix}
        \frac{n-1}{n} &-\frac{1}{\sqrt{2n}} &-\frac{1}{\sqrt{2n}} &\cdots &-\frac{1}{\sqrt{2n}}\\
        -\frac{1}{\sqrt{2n}} &\frac{1}{2} &0 &\cdots &0\\
        -\frac{1}{\sqrt{2n}} &0 &\frac{1}{2} &\cdots &0\\
        \cdots &\cdots &\cdots &\cdots &\cdots\\
        -\frac{1}{\sqrt{2n}} &0 &0 &\cdots &\frac{1}{2}
   \end{pmatrix}_{n \times n},
$$
\noindent and we easily have
$Sp(\mathcal{L}_A(S_n^2)) = \{0, \left(\frac{1}{2} \right)^{n-2}, \frac{3n-2}{2n} \}$.

Similarly, we have
$$\mathcal{L}_S(S_n^2) = \mathcal{L}_{11}(S_n^2) - \mathcal{L}_{12}(S_n^2) =
   \begin{pmatrix}
        \frac{n+1}{n} &-\frac{1}{\sqrt{2n}} &-\frac{1}{\sqrt{2n}} &\cdots &-\frac{1}{\sqrt{2n}}\\
        -\frac{1}{\sqrt{2n}} &\frac{3}{2} &0 &\cdots &0\\
        -\frac{1}{\sqrt{2n}} &0 &\frac{3}{2} &\cdots &0\\
        \cdots &\cdots &\cdots &\cdots &\cdots\\
        -\frac{1}{\sqrt{2n}} &0 &0 &\cdots &\frac{3}{2}
   \end{pmatrix}_{n \times n},
$$
\noindent and get
$Sp(\mathcal{L}_S(S_n^2)) = \{2, \left(\frac{3}{2} \right)^{n-2}, \frac{n+2}{2n} \}$.

Hence $Sp(\mathcal{L}(S_n^2)) = \{0, \left(\frac{1}{2} \right)^{n-2}, \frac{n+2}{2n},
 \frac{3n-2}{2n}, \left(\frac{3}{2} \right)^{n-2}, 2\}$ from Lemma 2.1,
and we immediately have the following result.

\noindent \textbf{Theorem 3.2.} Let $S_n^2 = S_n \times K_2$. Then

(1) $Kf^*(S_n^2) = \frac{48n^3 +25 n^2 -180n +116}{3n+6}$;

(2) $\lim\limits_{n \rightarrow +\infty} \frac{Kf^*(S_n^2)}{Gut(S_n^2)} = \frac{16}{33}$.

\noindent \textbf{Proof.} From Lemma 1.2 it is easily confirm that
\begin{align*}
 Kf^*(S_n^2) &= 2(3n-2) \left[ 2n-4 + \frac{2n}{n+2} + \frac{2n}{3n-2}+ \frac{2(n-2)}{3} +\frac{1}{2} \right] \\
 &= \frac{48n^3 +25 n^2 -180n +116}{3n+6}.
\end{align*}

Now, let $g_i =\sum_{j \in V(S_n^2)} d_i d_j d_{ij}$.
Obviously, if $i = 1, 1'$, then
$$g_i = n \cdot 2 \cdot 1 + n \cdot 2 \cdot 1 \cdot (n-1) + n \cdot 2 \cdot 2 \cdot (n-1)
      = 7n^2 -6n,$$
and otherwise
$$ g_i = 2 \cdot n \cdot 1 + 2 \cdot 2 \cdot 1 + 2 \cdot n \cdot 2 + 2 \cdot 2 \cdot 2 \cdot (n-2)
         + 2 \cdot 2 \cdot 3 \cdot (n-2) = 26n -36.$$

Hence
$$Gut(S_n^2) = \frac{1}{2}\sum\limits_{i \in V(S_n^2)} g_i
             = \frac{1}{2} \left[ 2(7n^2 -6n) + (26n -36)(2n-2) \right]
             = 33n^2 -68n +36,$$
and it follows that
$$\lim\limits_{n \rightarrow +\infty} \frac{Kf^*(S_n^2)}{Gut(S_n^2)}
   = \lim\limits_{n \rightarrow +\infty} \frac{48n^3 +25 n^2 -180n +116}{(3n+6)(33n^2 -68n +36)}
   =\frac{16}{33},$$
which completes the proof. $\blacksquare$

\section{Results for graphs in $\mathcal{S}_{n,r}^2$}
Let $S_{n,r}^2$ be any graph in $\mathcal{S}_{n,r}^2$, $1 \leq r \leq n-1$.
We will determine $Kf(S_{n,r}^2)$ and $\tau(S_{n,r}^2)$ in this section.

Let $d_i$ be the degree of vertex $i$ in $S_{n,r}^2$.
Then $d_i = n$ or $n-1$ if $i = 1, 1'$, and $d_i =1$ or $2$ otherwise.
We will compute $Sp(S_{n,r}^2)$ in the following two cases.

\textbf{Case 1.} Edge $11' \notin E(_{n,r}^2)$. Then
$$L_{11}(S_{n,r}^2) =\begin{pmatrix}
       n-1 &-1 &-1 &\cdots &-1 \\
       -1 &d_2 &0 &\cdots &0 \\
       -1 &0 &d_3 &\cdots &0 \\
       \cdots &\cdots &\cdots &\cdots &\cdots \\
       -1 &0 &0 &\cdots &d_n
    \end{pmatrix},~
  L_{12}(S_{n,r}^2) =\begin{pmatrix}
       0  &0 &0 &\cdots &0 \\
       0  &t_2 &0 &\cdots &0 \\
       0  &0   &t_3 &\cdots &0 \\
       \cdots &\cdots &\cdots &\cdots &\cdots \\
       0 &0 &0 &\cdots &t_n
    \end{pmatrix},
$$
where $t_i =0$ if $d_i =1$, and $t_i =1$ if $d_i =2$, $2 \leq i \leq n$.
Hence
$$L_A(S_{n,r}^2) = L_{11}(S_{n,r}^2) + L_{12}(S_{n,r}^2)
  =\begin{pmatrix}
       n-1 &-1 &-1 &\cdots &-1 \\
       -1 &1 &0 &\cdots &0 \\
       -1 &0 &1 &\cdots &0 \\
       \cdots &\cdots &\cdots &\cdots &\cdots \\
       -1 &0 &0 &\cdots &1
    \end{pmatrix}_{n \times n},
$$
and $Sp(L_A(S_{n,r}^2)) = \{0, 1^{n-2}, n \}$.

On the other hand,
$$L_S(S_{n,r}^2) = L_{11}(S_{n,r}^2) - L_{12}(S_{n,r}^2)
  =\begin{pmatrix}
       n-1 &-1 &-1 &\cdots &-1 \\
       -1 &d_2 -t_2 &0 &\cdots &0 \\
       -1 &0 &d_3 -t_3 &\cdots &0 \\
       \cdots &\cdots &\cdots &\cdots &\cdots \\
       -1 &0 &0 &\cdots &d_n -t_n
    \end{pmatrix},
$$
where $d_i -t_i =1$ if $d_i =1$, and $d_i -t_i =3$ if $d_i =2$, $2 \leq i \leq n$.
We will compute $Sp(L_S(S_{n,r}^2))$ in the following cases.

\textbf{Case 1.1.} $r=1$. Then $d_i -t_i =3$, $2 \leq i \leq n$,
and we easily have
$$ Sp(L_S(S_{n,r}^2))= \left\{ 3^{n-2}, \frac{n+2 +\sqrt{n^2 -4n +12}}{2}, \frac{n+2 -\sqrt{n^2 -4n +12}}{2}  \right\}.$$

\textbf{Case 1.2.} $r \geq 2$. By direct calculations, we have
$$\Phi( L_S(S_{n,r}^2), \lambda) =\left[\lambda^3 -(n+3)\lambda^2 +3n\lambda +2r -2n  \right]
   (\lambda -1)^{r-2} (\lambda -3)^{n-r-1}.$$

Let $\lambda_1, \lambda_2, \lambda_3$ be the three roots of $\lambda^3 -(n+3)\lambda^2 +3n\lambda +2r -2n = 0$.
Then $Sp(L_S(S_{n,r}^2)) =\{1^{r-2}, 3^{n-r-1}, \lambda_1, \lambda_2, \lambda_3 \}$,
and it holds that $\lambda_1 \lambda_2 \lambda_3 =2n -2r$
and
$$\frac{1}{\lambda_1} +\frac{1}{\lambda_2} +\frac{1}{\lambda_3}
 =\frac{\lambda_1 \lambda_2 +\lambda_1 \lambda_3 +\lambda_2 \lambda_3}{\lambda_1 \lambda_2 \lambda_3}
 =\frac{3n}{2n-2r}$$
from the Vieta's theorem.

\textbf{Case 2.} $11' \in E(_{n,r}^2)$. Then
$$L_{11}(S_{n,r}^2) =\begin{pmatrix}
       n  &-1 &-1 &\cdots &-1 \\
       -1 &d_2 &0 &\cdots &0 \\
       -1 &0 &d_3 &\cdots &0 \\
       \cdots &\cdots &\cdots &\cdots &\cdots \\
       -1 &0 &0 &\cdots &d_n
    \end{pmatrix},~
  L_{12}(S_{n,r}^2) =\begin{pmatrix}
       -1  &0 &0 &\cdots &0 \\
       0   &t_2 &0 &\cdots &0 \\
       0   &0   &t_3 &\cdots &0 \\
       \cdots &\cdots &\cdots &\cdots &\cdots \\
       0   &0 &0 &\cdots &t_n
    \end{pmatrix},
$$
where $t_i =0$ if $d_i =1$, and $t_i =-1$ if $d_i =2$, $2 \leq i \leq n$.
Hence
$$L_A(S_{n,r}^2) = L_{11}(S_{n,r}^2) + L_{12}(S_{n,r}^2)
  =\begin{pmatrix}
       n-1 &-1 &-1 &\cdots &-1 \\
       -1 &1 &0 &\cdots &0 \\
       -1 &0 &1 &\cdots &0 \\
       \cdots &\cdots &\cdots &\cdots &\cdots \\
       -1 &0 &0 &\cdots &1
    \end{pmatrix}_{n \times n},
$$
and $Sp(L_A(S_{n,r}^2)) = \{0, 1^{n-2}, n \}$.

On the other hand,
$$L_S(S_{n,r}^2) = L_{11}(S_{n,r}^2) - L_{12}(S_{n,r}^2)
  =\begin{pmatrix}
       n+1 &-1 &-1 &\cdots &-1 \\
       -1 &d_2 -t_2 &0 &\cdots &0 \\
       -1 &0 &d_3 -t_3 &\cdots &0 \\
       \cdots &\cdots &\cdots &\cdots &\cdots \\
       -1 &0 &0 &\cdots &d_n -t_n
    \end{pmatrix},
$$
where $d_i -t_i =1$ if $d_i =1$, and $d_i -t_i =3$ if $d_i =2$, $2 \leq i \leq n$.
By direct calculations, we have
$$\Phi( L_S(S_{n,r}^2), \lambda) =\left[\lambda^3 -(n+5)\lambda^2 +(3n+8)\lambda +2r -2n -4 \right]
   (\lambda -1)^{r-1} (\lambda -3)^{n-r-2}.$$

Let $\lambda_1, \lambda_2, \lambda_3$ be the three roots of $\lambda^3 -(n+5)\lambda^2 +(3n+8)\lambda +2r -2n -4 = 0$.
Then $Sp(L_S(S_{n,r}^2)) =\{1^{r-1}, 3^{n-r-2}, \lambda_1, \lambda_2, \lambda_3 \}$,
and it holds that $\lambda_1 \lambda_2 \lambda_3 =2n -2r +4$ and
$$\frac{1}{\lambda_1} +\frac{1}{\lambda_2} +\frac{1}{\lambda_3}
 =\frac{\lambda_1 \lambda_2 +\lambda_1 \lambda_3 +\lambda_2 \lambda_3}{\lambda_1 \lambda_2 \lambda_3}
 =\frac{3n+8}{2n-2r+4}$$ from the Vieta's theorem.

Now, we are able to give the main result of this section.

\noindent \textbf{Theorem 4.1.} If $S_{n,r}^2 \in \mathcal{S}_{n,r}^2$, $0 \leq r \leq n-1$, then

(1) $Kf(S_{n,r}^2) = \begin{cases}
  \frac{8n^3 -(4r+17)n^2 -(4r^2 -26r -6)n -6r}{3(n-r)}, &\text{if}~11' \notin E(S_{n,r}^2) \\
  \frac{8n^3 -(4r-3)n^2 -(4r^2 -30r +14)n +12 -6r}{3(n-r-2)}, &\text{if}~11' \in E(S_{n,r}^2)
\end{cases} $;

(2) $\tau(S_{n,r}^2) = \begin{cases}
  (n-r) \cdot 3^{n-r-1}, &\text{if}~11' \notin E(S_{n,r}^2) \\
  (n-r+2) \cdot 3^{n-r+2}, &\text{if}~11' \in E(S_{n,r}^2)
\end{cases} $;

(3) $\lim\limits_{n \rightarrow +\infty} \frac{Kf(S_{n,r}^2)}{W(S_{n,r}^2)} = \frac{8}{15}$.

\noindent \textbf{Proof.} If $r=0$, then $S_{n,r}^2 \cong S_{n}^2$, and the conclusion holds from Theorem 3.1.
Hence assume $r \geq 1$. We distinguish the following two cases.

\textbf{Case 1.} Edge $11' \notin E(_{n,r}^2)$.

\textbf{Case 1.1.} $r=1$. Then
$$ Sp(L(S_{n,r}^2)) =\{0, 1^{n-2}, n, 3^{n-2},
  \frac{n+2 -\sqrt{n^2 -4n +12}}{2}, \frac{n+2 +\sqrt{n^2 -4n +12}}{2} \},$$
From Lemma 1.1 we have
\begin{align*}
 Kf(S_{n,r}^2) &= 2n \left[n-2 +\frac{1}{n} +\frac{n-2}{3}
    +\frac{2}{n+2 -\sqrt{n^2 -4n +12}} +\frac{2}{n+2 +\sqrt{n^2 -4n +12}} \right] \\
  &= \frac{8n^3 -21n^2 +28n -6}{3(n-1)}\\
  &= \frac{8n^3 -(4r+17)n^2 -(4r^2 -26r -6)n -6r}{3(n-r)}.
\end{align*}

Then from Lemma 2.2, we have
\begin{align*}
 \tau(S_{n,r}^2) &= \frac{1}{2n} \left[n \cdot 3^{n-2}
    \cdot \frac{n+2 -\sqrt{n^2 -4n +12}}{2} \cdot \frac{n+2 +\sqrt{n^2 -4n +12}}{2} \right] \\
  &= (n-1) \cdot 3^{n-2}\\
  &= (n-r) \cdot 3^{n-r-1}.
\end{align*}

\textbf{Case 1.2.} $r \geq 2$. Then $ Sp(L(S_{n,r}^2)) =\{0, 1^{n+r-4}, n, 3^{n-r-1}, \lambda_1, \lambda_2, \lambda_3 \}$, where $\lambda_1 \lambda_2 \lambda_3 = 2n-2r$
and $1/\lambda_1 + 1/\lambda_2 + 1/\lambda_3 = 3n/(2n-2r)$.
From Lemma 1.1 we have
\begin{align*}
 Kf(S_{n,r}^2) &= 2n \left[n+r-4 +\frac{1}{n} +\frac{n-r-1}{3} + \frac{3n}{2n-2r} \right] \\
  &= \frac{8n^3 -(4r+17)n^2 -(4r^2 -26r -6)n -6r}{3(n-r)}.
\end{align*}

Then from Lemma 2.2, we have
$$  \tau(S_{n,r}^2) = \frac{n \cdot 3^{n-r-1} \cdot \lambda_1 \cdot \lambda_2 \cdot \lambda_3}{2n}
  = \frac{n \cdot 3^{n-r-1} \cdot (2n-2r)}{2n}
  = (n-r) \cdot 3^{n-r-1}.
$$

\textbf{Case 2.} Edge $11' \in E(_{n,r}^2)$. Then $ Sp(L(S_{n,r}^2)) =\{0, 1^{n+r-3}, n, 3^{n-r-2}, \lambda_1, \lambda_2, \lambda_3 \}$, where $\lambda_1 \lambda_2 \lambda_3 = 2n-2r+4$
and $1/\lambda_1 + 1/\lambda_2 + 1/\lambda_3 = (3n+8)/(2n-2r+4)$.
From Lemma 1.1 we have
\begin{align*}
 Kf(S_{n,r}^2) &= 2n \left[n+r-3 +\frac{1}{n} +\frac{n-r-2}{3} + \frac{3n+8}{2n-2r+4} \right] \\
  &= \frac{8n^3 -(4r-3)n^2 -(4r^2 -30r +14)n +12 -6r}{3(n-r+2)}.
\end{align*}

Then from Lemma 2.2, we have
$$  \tau(S_{n,r}^2) = \frac{n \cdot 3^{n-r-2} \cdot \lambda_1 \cdot \lambda_2 \cdot \lambda_3}{2n}
  = \frac{n \cdot 3^{n-r-2} \cdot (2n-2r+4)}{2n}
  = (n-r+2) \cdot 3^{n-r-2}.
$$

Finally, it is straightforward to have $W(S_{n,r}^2) = W(S_{n}^2) +r = 5n^2 -8n +r+4$.
Hence, in both cases, it holds that
$$\lim\limits_{n \rightarrow +\infty} \frac{Kf(S_{n,r}^2)}{W(S_{n,r}^2)} = \frac{8}{15}.~\blacksquare$$

\end{document}